\documentclass[review]{elsarticle}

\usepackage{hyperref}

\usepackage{graphics,psfrag,fullpage}
\usepackage{graphicx, amssymb, amsmath}
\usepackage{color}
\usepackage{url}
\usepackage{amsmath}
\numberwithin{equation}{section}

\newtheorem{theorem}{Theorem}[section]

\newtheorem{lemma}{Lemma}[section]
\newtheorem{remark}{Remark}[section]
\newtheorem{examp}{Example}[section]

\journal{Journal of \LaTeX\ Templates}

\begin{document}

\begin{frontmatter}

\title{Stability criteria of linear delay differential systems based on fundamental matrix}

\author{Guang-Da Hu*}
\address{Shanghai Customs University, Shanghai, 201204, China\\
and\\
Department of Mathematics\\
Shanghai University, Shanghai, 200444, China}
\ead{ghu@hit.edu.cn}


\cortext[mycorrespondingauthor]{Guang-Da Hu}

\begin{abstract}
We investigate stability of linear delay differential systems. Stability criteria of the systems
are derived based on integrals of the fundamental matrix. They are necessary and sufficient conditions
for delay-dependent stability of the systems. Numerical examples are given to illustrate the main results.
\end{abstract}


\end{frontmatter}

\section{Introduction}
We investigate delay-dependent stability of linear delay differential systems described by
\begin{equation}
\dot{x}(t)=A_{0}x(t)+\sum_{j=1}^{m}A_{j}x(t-\tau_{j}),
\label{eq:dde}
\end{equation}
where constant matrices $A_{0}, A_{j}\in {\mathcal{R}}^{n\times n}$ and time delays $\tau_{j}>0$ for $j=1,\ldots, m.$
\par
Delay differential systems arise in different applied problems \cite{Kol}. Delays in engineering systems can have several causes \cite{Kol,Gu},
for example: technological lag, signal transmission and information delay, time of mixing reactants. Delay differential systems may also be applied to describe economic dynamics \cite{Kel}.
\par
The stability problem of linear delay differential systems is a challenging research topic which has received much attention
in the references \cite{Kol,Gu,Fri,Khari,Kim}. The large number of approaches
in the references yielding delay-dependent stability conditions can be classified into two main branches.
One branch is based on the construction of a Lyapunov-Krasovskii functional,
the other is based on the discussion of characteristic roots. In \cite{Gu,Fri,Card}, for the stability of system (\ref{eq:dde}),
a Lyapunov-Krasovskii functional is constructed and sufficient but conservative delay-dependent stability conditions are
formulated as linear matrix inequalities (LMIs). An advantage of this approach is that conditions can be formulated with LMIs which allow the automatic application in engineering software.
The disadvantage is that the results may be conservative. The approaches based on the characteristic roots typically yield exact results
but normally complicated computations are needed \cite{Gu,Sch}.
\par
Recently, based on the fundamental matrix, stability of system (\ref{eq:dde}) has been
investigated in \cite{Gy,Hu1,Hu2,Sab}, respectively. Based on integrals of the fundamental matrix,
necessary conditions of stability of system (\ref{eq:dde}) are given in \cite{Gy,Sab}. In \cite{Hu1}, a stability criterion is presented via
the integral of the square of $F$-norm for the fundamental matrix. Further, the stability criterion is applied
to design stabilizing controllers of system (\ref{eq:dde}) in \cite{Hu1,Hu2}.
\par
Along the line of \cite{Gy,Hu1,Hu2,Sab}, the main contributions are summarized as follows:
\begin{description}
\item{1.}
Based on the fundamental matrix, stability criteria are derived which are necessary and sufficient conditions
for delay-dependent stability of (\ref{eq:dde}).
\item{2.}
We prove that the necessary condition of stability of (\ref{eq:dde}) in \cite{Gy,Sab} is also sufficient.
\end{description}
\par
Throughout this paper, the absolute value of a real number $w$ is denoted by $|w|$.
$||W||$ and $||W||_{F}$ stand for the norm and the Frobenius norm of a matrix $W,$
respectively.
\section{Preliminaries}
Now we review the concept of the fundamental matrix of system (\ref{eq:dde}).
The fundamental matrix of system (\ref{eq:dde}), $X(t)\in {\mathcal{R}}^{n\times n}$
satisfies system (\ref{eq:dde}), {\it i. e.,}
\begin{equation}
\dot{X}(t)=A_{0}X(t)+\sum_{j=1}^{m}A_{j}X(t-\tau_{j}) \quad \mbox{for} \quad t\geq 0,
\label{eq:state1}
\end{equation}
under the condition
\begin{equation}
X(t)=0 \quad \mbox{for} \quad t<0 \quad \mbox{and}\quad X(0)=I.
\label{eq:state2}
\end{equation}
We know that system (\ref{eq:dde}) is asymptotically stable if and only if
 \[
||X(t)||\leq \gamma\exp(-\alpha t), \quad \mbox{for}\quad t\in [0, \infty),
\]
where $\alpha$ and $\gamma>0.$ Thus stability of system (\ref{eq:dde}) is completely determined by the properties of the fundamental matrix \cite{Kol,Gu,Fri,Kim}.
\par
Based on the fundamental matrix, the following result has been obtained in \cite{Gy,Sab}.
\begin{lemma}(\cite{Gy,Sab})
If system (\ref{eq:dde}) is asymptotically stable, then the matrix
$(A_{0}+\sum_{j=1}^{m}A_{j})$ is invertible and
\[
\int_{0}^{\infty} X(\sigma)d \sigma=-(A_{0}+\sum_{j=1}^{m}A_{j})^{-1}.
\]
\label{lem:Nec}
\end{lemma}
\par
In order to numerically solve the functional matrix $X(t),$ it is necessary to discretize it.
Assume the numerical solution gives a sequence of approximated values\\
$\{
X_{1}, X_{2}, \ldots, X_{N}\}$ of $\{X(t_{1}), X(t_{2}), \ldots, X(t_{N})\}$ of Eq. (\ref{eq:state1}) on certain equidistant step-values $\{t_{n}=nh\}$ with the step-size $h=T/N,$ where $N$ is a positive integer. Let $l_{i}=[\tau_{i}h^{-1}],\delta_{i}=l_{i}-\tau_{i}h^{-1},0\leq \delta_{i}<1,$ for $i=1, \ldots, m,$ where $[a]$ stand for the smallest integer that is greater than or equal to $a\in \mathcal{R}$.
%
Using Euler scheme combining with the linear interpolation for Eqs.(\ref{eq:state1}) and
(\ref{eq:state2}), for $ n=0,1,\ldots,N-1,$ we have the following scheme
\begin{equation}
X_{n+1}=X_{n}+h(A_{0}X_{n}+\sum_{j=1}^{m}A_{j}X_{n-l_{j}+\delta_{j}}),
\label{eq:state1-d}
\end{equation}
\begin{equation}
X_{n-l_{j}+\delta_{j}}=0 \quad \mbox{for} \quad n-l_{j}+\delta_{j}<0,
\label{eq:state2-d}
\end{equation}
\begin{equation}
X_{n-l_{j}+\delta_{j}}=I \quad \mbox{for} \quad n-l_{j}+\delta_{j}=0,
\label{eq:state3-d}
\end{equation}
and
\begin{equation}
X_{n-l_{j}+\delta_{j}}=(1-\delta_{j})X_{n-l_{j}}+\delta_{j}X_{n-l_{j}+1}
\quad \mbox{for} \quad n-l_{j}+\delta_{j}>0.
\label{eq:state4-d}
\end{equation}
\par
Finally, we review the following result which is known as Barbalat's lemma.
\begin{lemma}(\cite{Kha})
\label{lem:BL}
Let $f(t)$ be uniformly continuous function on $[0, \infty).$ Suppose that
\[
\int_{0}^{\infty}f(\sigma)d \sigma
\]
exists and is finite. Then
\[
f(t)\rightarrow 0\quad \mbox{for}\quad t\rightarrow \infty.
\]
\label{lem:BL}
\end{lemma}

\section{Stability Criteria}
In this section, based on the fundamental matrix $X(t),$ we present main results which
are necessary and sufficient conditions for delay-dependent stability of system (\ref{eq:dde}).
\par
\begin{theorem}
\label{th:P1}
System (\ref{eq:dde}) is asymptotically stable if and only if
\begin{equation}
A_{0}+\sum_{j=1}^{m}A_{j}\quad \mbox{is invertible and}\quad
\int_{0}^{\infty} X(\sigma)d \sigma=-(A_{0}+\sum_{j=1}^{m}A_{j})^{-1}.
\label{eq:sqe}
\end{equation}
\end{theorem}
{\it Proof.}
Suppose that condition (\ref{eq:sqe}) holds.
\par
Integrating both sides of system (\ref{eq:state1}) from $0$ to $t,$ we have that
 \begin{equation}
\int_{0}^{t}\dot{X}(\sigma)d\sigma
=A_{0}\int_{0}^{t}X(\sigma)d\sigma+\sum_{j=1}^{m}A_{j}\int_{0}^{t}
X(\sigma-\tau_{j})d \sigma.
\label{eq:insq11}
\end{equation}
For $X(\sigma)=0$ for $\sigma<0,$ the right side of Eq. (\ref{eq:insq11})
\begin{eqnarray}
\begin{split}
&A_{0}\int_{0}^{t}X(\sigma)d\sigma+\sum_{j=1}^{m}A_{j}\int_{0}^{t}X(\sigma-\tau_{j})d \sigma\\
&=A_{0}\int_{0}^{t}X(\sigma)d\sigma+\sum_{j=1}^{m}A_{j}\int_{-\tau_{j}}^{t-\tau_{j}}X(\sigma)d \sigma\\
&=A_{0}\int_{0}^{t}X(\sigma)d\sigma+\sum_{j=1}^{m}A_{j}\int_{-\tau_{j}}^{0}X(\sigma)d\sigma
+\sum_{j=1}^{m}A_{j}\int_{0}^{t-\tau_{j}}X(\sigma)d \sigma.
\end{split}
\label{eq:insq12}
\end{eqnarray}
Since $X(\sigma)=0$ for $\sigma<0$, in the right side of Eq. (\ref{eq:insq12}), $\sum_{j=1}^{m}A_{j}\int_{-\tau_{j}}^{0}X(\sigma)d\sigma=0.$
Let $t \rightarrow \infty$ in Eq. (\ref{eq:insq12}), we obtain that
\begin{eqnarray*}
&&A_{0}\int_{0}^{\infty}X(\sigma)d\sigma+\sum_{j=1}^{m}A_{j}\int_{0}^{\infty}X(\sigma-\tau_{j})d \sigma\\
&&=A_{0}\int_{0}^{\infty}X(\sigma)d\sigma+\sum_{j=1}^{m}A_{j}\int_{0}^{\infty}X(\sigma)d \sigma\\
&&=(A_{0}+\sum_{j=1}^{m}A_{j})\int_{0}^{\infty}X(\sigma)d \sigma.
\end{eqnarray*}
According to (\ref{eq:sqe}), we have
\begin{equation}
\begin{split}
&A_{0}\int_{0}^{\infty}X(\sigma)d\sigma+\sum_{j=1}^{m}A_{j}\int_{0}^{\infty}X(\sigma-\tau_{j})d \sigma\\
&=(A_{0}+\sum_{j=1}^{m}A_{j})\int_{0}^{\infty}X(\sigma)d \sigma=-I.
\end{split}
\label{eq:insq22}
\end{equation}
It means the generalized integral on the left side of Eq. (\ref{eq:insq11}) exists and
\begin{equation}
 \int_{0}^{\infty}\dot{X}(\sigma)d\sigma=\lim_{t \rightarrow \infty}X(t)-X(0).
 \label{eq:insq3}
\end{equation}
By means of Eqs. (\ref{eq:insq22}) and (\ref{eq:insq3}),
\[
\lim_{t \rightarrow \infty}X(t)-X(0)=-I.
\]
Since $X(0)=I,$
\[
\lim_{t \rightarrow \infty}X(t)=0,
\]
{\it i. e.,} system (\ref{eq:dde}) is asymptotically stable.
\par
Conversely, suppose that system (\ref{eq:dde}) is asymptotically stable.
We have that $\lim_{t \rightarrow \infty}X(t)=0.$
\par
Let $t \rightarrow \infty$ in the left side of Eq. (\ref{eq:insq11}), for $X(0)=I,$ we obtain that
\[
\int_{0}^{\infty}\dot{X}(\sigma)d\sigma=\lim_{t \rightarrow \infty}X(t)-X(0)=-I.
\]
Using the above equation and $X(\sigma)=0$ for $\sigma<0$, as $t \rightarrow \infty,$ the generalized integral
in the right side Eq. (\ref{eq:insq11}) exists and
\[
-I=(A_{0}+\sum_{j=1}^{m}A_{j})\int_{0}^{\infty}X(\sigma)d \sigma.
\]
\par
The characteristic equation of system (\ref{eq:dde}) is
\[
P(s)=\det[sI-A_{0}-\sum_{j=1}^{m}A_{j}\exp(-\tau_{j}s)]=0.
\]
Since system (\ref{eq:dde}) is asymptotically stable, $s=0$ is not a root of the characteristic equation.
Thus
\[
P(0)\neq 0\quad \Rightarrow \quad \det[A_{0}+\sum_{j=1}^{m}A_{j}]\neq 0.
\]
It means that the matrix $(A_{0}+\sum_{j=1}^{m}A_{j})$ is invertible.
\par
The proof is completed. $\square$
\par
\begin{remark}
The necessity of condition (\ref{eq:sqe}) for stability of system (\ref{eq:dde})
has been given in \cite{Gy,Sab}. See Lemma \ref{lem:Nec}. Theorem \ref{th:P1} shows that
 condition (\ref{eq:sqe}) is also sufficient for stability of system (\ref{eq:dde}).
\end{remark}
\par
\begin{remark}
For linear systems without delays
\begin{equation}
\dot{x}(t)=Ax(t),
\label{eq:ode}
\end{equation}
where $X(t)=e^{At}$ is the fundamental matrix of system (\ref{eq:ode}).
Condition (\ref{eq:sqe}) in Theorem \ref{th:P1} reduces to
\[
A\quad \mbox{is invertible and}\quad
\int_{0}^{\infty} X(\sigma)d \sigma=-A^{-1}.
\]
\end{remark}
\par
\begin{theorem}
\label{th:P2}
System (\ref{eq:dde}) is asymptotically stable if and only if
\begin{equation}
\int_{0}^{\infty} {||X(\sigma)||} d \sigma \quad \mbox{exists and is finite}.
\label{eq:int-norm}
\end{equation}
\end{theorem}
{\it Proof.}
Suppose that condition (\ref{eq:int-norm}) holds.
\par
According to Eq. (\ref{eq:insq11}) and $X(\sigma)=0$ for $\sigma<0$, we obtain that
\begin{eqnarray*}
&&||\int_{0}^{t}\dot{X}(\sigma)d\sigma||\\
&&\leq || A_{0}\int_{0}^{t}X(\sigma)d\sigma||+\sum_{j=1}^{m}||A_{j}\int_{0}^{t-\tau_{j}}X(\sigma)d \sigma||\\
&&\leq || A_{0}||\int_{0}^{t}||X(\sigma)||d\sigma||+\sum_{j=1}^{m}||A_{j}||\int_{0}^{t-\tau_{j}}||X(\sigma)||d \sigma\\
&&\leq || A_{0}||\int_{0}^{t}||X(\sigma)||d\sigma||+\sum_{j=1}^{m}||A_{j}||\int_{0}^{t}||X(\sigma)||d \sigma.
\end{eqnarray*}
That is
\begin{equation}
||\int_{0}^{t}\dot{X}(\sigma)d\sigma||\leq (|| A_{0}||+\sum_{j=1}^{m}||A_{j}||)\int_{0}^{t}||X(\sigma)||d \sigma.
\label{eq:int-norm2}
\end{equation}
Since $\int_{0}^{\infty} {||X(\sigma)||} d \sigma $ is bounded, from (\ref{eq:int-norm2}) we know that
 $||\int_{0}^{t}\dot{X}(\sigma)d\sigma||$ is bounded. Thus
 \[
 ||\int_{0}^{t}\dot{X}(\sigma)d\sigma||=||X(t)-X(0)||=||X(t)-I||
 \]
 is bounded. The above inequality means that $||X(t)||$ is bounded for $t\in [0, \infty).$
\par
By means of Eq. (\ref{eq:state1}),
\begin{equation}
||\dot{X}(t)||\leq ||A_{0}|| ||X(t)||+\sum_{j=1}^{m}||A_{j}|| ||X(t-\tau_{j})||\quad \mbox{for} \quad t\geq 0.
\label{eq:int-n1}
\end{equation}
According to $||X(t)||$ is bounded for $t\in [0, \infty)$ and (\ref{eq:int-n1}),
we know that $||\dot{X}(t)||$ is bounded for $t\in [0, \infty).$ It implies that
$\dot{X}_{ij}(t)$ is bounded for $t\in [0, \infty),$ where $X(t)=[X_{ij}(t)].$ Thus $X_{ij}(t)$ is uniformly continuous for $t\in [0, \infty).$
Since $\int_{0}^{\infty} {||X(\sigma)||} d \sigma$ exists and is finite, $\int_{0}^{\infty} X(\sigma) d \sigma $ exists and is finite.
By means of Lemma \ref{lem:BL}, we obtain that
\[
X_{ij}(t) \rightarrow 0\quad \mbox{as} \quad t\rightarrow \infty.
\]
Thus
\[
||X(t)||\rightarrow 0\quad \mbox{as} \quad t\rightarrow \infty,
\]
%
{\it i. e.,} system (\ref{eq:dde}) is asymptotically stable.
\par
Conversely, suppose that system (\ref{eq:dde}) is asymptotically stable. We know that
 \[
||X(t)||\leq \gamma\exp(-\alpha t), \quad \mbox{for}\quad t\in [0, \infty),
\]
where $\alpha$ and $\gamma$ are positive constants. It is obvious that $\int_{0}^{\infty} {||X(\sigma)||} d \sigma $ exists and is finite.
\par
The proof is completed. $\square$
\par
\begin{remark}
For system (\ref{eq:ode}), we have the following result:
system (\ref{eq:ode}) is asymptotically stable if and only if
\[
\int_{0}^{\infty} {||X(\sigma)||} d \sigma \quad \mbox{exists and is finite}.
\]
See \cite{Rugh}. Thus Theorem \ref{th:P2} is an extension of the result to linear delay systems.
\end{remark}
\par
\begin{theorem}
\label{th:P3}
System (\ref{eq:dde}) is asymptotically stable if and only if
\begin{equation}
\int_{0}^{\infty} {||X(\sigma)||_{F}}^{2} d \sigma \quad \mbox{exists and is finite.}
\label{eq:int-normf}
\end{equation}
\end{theorem}
{\it Proof.}
Suppose that condition (\ref{eq:int-normf}) holds.
\par
Let matrices $A_{0}=[a_{ij}],$ $A_{l}=[a_{ij}^{(l)}]$ and $X(t)=[X_{ij}(t)].$
System (\ref{eq:state1}) is rewritten as
\[
\dot{X}_{ij}(t)=\sum_{k=1}^{n}[a_{ik}X_{kj}(t)+\sum_{l=1}^{m}{a_{ik}}^{(l)}X_{kj}(t-\tau_{l})]
\quad \mbox{for} \quad i,j=1, \ldots,n,\quad t\geq 0.
\]
Multiplying the two sides of the above equation by $X_{ij}(t),$ for $i,j=1, \ldots,n, t\geq 0,$
\begin{equation}
X_{ij}(t)\dot{X}_{ij}(t)=\sum_{k=1}^{n}[a_{ik}X_{ij}(t)X_{kj}(t)+\sum_{l=1}^{m}{a_{ik}}^{(l)}X_{ij}(t)X_{kj}(t-\tau_{l})].
\label{eq:sq1}
\end{equation}
Integrating both sides of system (\ref{eq:sq1}) from $0$ to $t,$ we have that for $i,j=1, \ldots,n, t\geq 0,$
\begin{equation}
\int_{0}^{t}X_{ij}(\sigma)\dot{X}_{ij}(\sigma)d\sigma=\sum_{k=1}^{n}[a_{ik}\int_{0}^{t}X_{ij}(\sigma)X_{kj}(\sigma)d\sigma+\sum_{l=1}^{m}
{a_{ik}}^{(l)}\int_{0}^{t}X_{ij}(\sigma)X_{kj}(\sigma-\tau_{l})d\sigma].
\label{eq:insq1}
\end{equation}
Since
\begin{eqnarray*}
|X_{ij}(\sigma)X_{kj}(\sigma)|\leq  {X_{ij}(\sigma)}^{2}+{X_{kj}(\sigma)}^{2}\quad \mbox{and} \quad
|X_{ij}(\sigma)X_{kj}(\sigma-\tau_{l})|\leq {X_{ij}(\sigma)}^{2}+{X_{kj}(\sigma-\tau_{l})}^{2},
\end{eqnarray*}
we have that
 \begin{equation}
\int_{0}^{t}|X_{ij}(\sigma)X_{kj}(\sigma)|d\sigma \leq \int_{0}^{t}{X_{ij}(\sigma)}^{2}d\sigma+\int_{0}^{t}{X_{kj}(\sigma)}^{2}d\sigma
\label{eq:insqa}
\end{equation}
and
\begin{eqnarray*}
&&\int_{0}^{t}|X_{ij}(\sigma)X_{kj}(\sigma-\tau_{l})|d\sigma\leq \int_{0}^{t}[{X_{ij}(\sigma)}^{2}+{X_{kj}(\sigma-\tau_{l})}^{2}]d\sigma\\
&&=\int_{0}^{t}{X_{ij}(\sigma)}^{2}d\sigma + \int_{0}^{t}{X_{kj}(\sigma-\tau_{l})}^{2}d\sigma \\
&&=\int_{0}^{t}{X_{ij}(\sigma)}^{2}d\sigma + \int_{0}^{t-\tau_{l}}{X_{kj}(\sigma)}^{2}d\sigma +
\int_{-\tau_{l}}^{0}{X_{kj}(\sigma)}^{2}d\sigma\\
&&\leq \int_{0}^{t}{X_{ij}(\sigma)}^{2}d\sigma + \int_{0}^{t}{X_{kj}(\sigma)}^{2}d\sigma +
\int_{-\tau_{l}}^{0}{X_{kj}(\sigma)}^{2}d\sigma.
\end{eqnarray*}
Since $X_{kj}(\sigma)=0$ for $\sigma<0,$ we have that
 \begin{equation}
 \int_{0}^{t}|X_{ij}(\sigma)X_{kj}(\sigma-\tau_{l})|d\sigma
 \leq \int_{0}^{t}{X_{ij}(\sigma)}^{2}d\sigma +\int_{0}^{t}{X_{kj}(\sigma)}^{2}d\sigma.
\label{eq:insqb}
\end{equation}
Using (\ref{eq:insq1}), (\ref{eq:insqa}) and (\ref{eq:insqb}), we have that
\begin{eqnarray*}
&&|\int_{0}^{t}d {X_{ij}}^{2}(\sigma)|=2|\int_{0}^{t}X_{ij}(\sigma)\dot{X}_{ij}(\sigma)d\sigma|\\
&&\leq 2 \sum_{k=1}^{n}[|a_{ik}|\int_{0}^{t}|X_{ij}(\sigma)X_{kj}(\sigma)|d\sigma+
\sum_{l=1}^{m}{|a_{ik}}^{(l)}|\int_{0}^{t}|X_{ij}(\sigma)X_{kj}(\sigma-\tau_{l})|d\sigma]\\
&& \leq 2 \sum_{k=1}^{n}\{|a_{ik}|[\int_{0}^{t}{X_{ij}(\sigma)}^{2}d\sigma+\int_{0}^{t} {X_{kj}(\sigma)}^{2}d\sigma]\\
&& +\sum_{l=1}^{m}{|a_{ik}}^{(l)}|[\int_{0}^{t}{X_{ij}(\sigma)}^{2}d\sigma+\int_{0}^{t}{X_{kj}(\sigma)}^{2}d\sigma]\}.
\end{eqnarray*}
From $\int_{0}^{\infty} {||X(\sigma)||_{F}}^{2} d \sigma$ is bounded, we know that
$\int_{0}^{\infty} {X_{ij}(\sigma)}^{2} d \sigma$ is bounded for $i,j=1, \ldots,n.$ Thus
$|\int_{0}^{\infty}d {X_{ij}}^{2}(\sigma)|$ is bounded. Since
\[
|\int_{0}^{t}d {X_{ij}}^{2}(\sigma)|=|{X_{ij}}^{2}(t)-{X_{ij}}^{2}(0)|,
\]
we obtain that
\[
|\int_{0}^{\infty}d {X_{ij}}^{2}(\sigma)|=\lim_{t \rightarrow \infty}|{X_{ij}}^{2}(t)-{X_{ij}}^{2}(0)|
\]
is bounded, {\it i. e.,} ${X_{ij}}^{2}(t)$ is bounded for $i,j=1, \ldots,n,t\in [0, \infty).$
It implies that  ${X_{ij}}(t)$ is bounded for $i,j=1, \ldots,n,t\in [0, \infty),$ and $||X(t)||$ is bounded.
\par
According to $||X(t)||$ is bounded for $t\in [0, \infty)$ and inequality (\ref{eq:int-n1}),
we know that $||\dot{X}(t)||$ is bounded for $t\in [0, \infty).$ It implies that
$\dot{X}_{ij}(t)$ is bounded for $i,j=1, \ldots,n,t\in [0, \infty).$ Since
${X_{ij}}(t)$ and $\dot{X}_{ij}(t)$ are bounded $i,j=1, \ldots,n,t\in [0, \infty),$
let $|{X_{ij}}(t)|<c_{1}$ and $|\dot{X}_{ij}(t)|<c_{2},$ where $c_{1}$ and $c_{2}$ are
positive constants. By the mean value theorem, for any $t_{1}$ and $t_{2}\in [0, \infty)$
with $t_{2}>t_{1}$, there is a $\xi$ with $t_{1}<\xi<t_{2}$ such that
\[
{X_{ij}}^{2}(t_{2})-{X_{ij}}^{2}(t_{1})=2X_{ij}(\xi)\dot{X}_{ij}(\xi)(t_{2}-t_{1}),\quad \mbox{for} \quad i,j=1, \ldots,n,t\in [0, \infty).
\]
Thus
\[
|{X_{ij}}^{2}(t_{2})-{X_{ij}}^{2}(t_{1})|=2|X_{ij}(\xi)\dot{X}_{ij}(\xi)| |t_{2}-t_{1}|< 2c_{1}c_{2} |t_{2}-t_{1}|,
\quad \mbox{for} \quad i,j=1, \ldots,n,t\in [0, \infty).
\]
It means that ${X_{ij}(t)}^{2}$ is uniformly continuous for $i,j=1, \ldots,n,t\in [0, \infty).$
Thus ${||X(t)||_{F}}^{2}$ is uniformly continuous for $i,j=1, \ldots,n,t\in [0, \infty).$ By means of Lemma \ref{lem:BL}, we obtain that
\[
{||X(t)||_{F}}^{2}\rightarrow 0\quad \mbox{as} \quad t\rightarrow \infty
\]
and
\[
{||X(t)||_{F}}\rightarrow 0\quad \mbox{as} \quad t\rightarrow \infty.
\]
System (\ref{eq:dde}) is asymptotically stable.
\par
Conversely, suppose that system (\ref{eq:dde}) is asymptotically stable. We know that
 \[
||X(t)||_{F}\leq \gamma \exp(-\alpha t), \quad \mbox{for}\quad t\in [0, \infty),
\]
where $\alpha$ and $\gamma$ are positive constants.
It is obvious that $\int_{0}^{\infty} {||X(\sigma)||_{F}}^{2} d \sigma$ exists and is finite.
\par
The proof is completed. $\square$
\par
\begin{remark}
A proof of Theorem \ref{th:P3} has been given in \cite{Hu1} which relies on
the uniform continuity of the fundamental matrix. In the present paper,
the proof is elementary which needs only knowledge of
classical mathematical analysis. Theorem \ref{th:P3} has been applied
to design stabilizing controllers of system (\ref{eq:dde}) in \cite{Hu1,Hu2,Hu3,Hu4}.
\end{remark}
\section{Numerical Examples}
In this section, numerical examples are given to illustrate the main results.
\par
For system (\ref{eq:dde}), using the scheme, Eqs. (\ref{eq:state1-d}) to (\ref{eq:state4-d}), we
can obtain a numerical solution $X_{n}$ of the fundamental matrix $X(t).$  By means of
the simplest numerical integration, we have approximated values
\begin{equation}
\sum_{n=1}^{N}X_{n}h \quad \mbox{for} \quad \int_{0}^{\infty} X(\sigma)d \sigma,
\label{eq:INT-1}
\end{equation}
\begin{equation}
\sum_{n=1}^{N}||X_{n}||_{F} h \quad \mbox{for} \quad \int_{0}^{\infty} {||X(\sigma)||_{F}} d \sigma,
\label{eq:INT-2}
\end{equation}
and
\begin{equation}
\sum_{n=1}^{N}{||X_{n}||_{F}}^{2} h \quad \mbox{for} \quad \int_{0}^{\infty} {||X(\sigma)||_{F}}^{2} d \sigma,
\label{eq:INT-3}
\end{equation}
respectively. Eqs. (\ref{eq:INT-1}), (\ref{eq:INT-2}) and (\ref{eq:INT-3}) are used to verify
Theorems \ref{th:P1}, \ref{th:P2}  and \ref{th:P3}, respectively.
\par
\begin{examp}
For system (\ref{eq:dde}), let $m=2.$ The system matrices are
$$A_{0}=\left[\begin{array}{ccc}
-4.2583&-0.7236&-11.8071\\
-5.6881&-2.7135&-6.8457\\
1.2493&-1.0037&-6.3093
\end{array}\right],
A_{1}=\left[\begin{array}{ccc}
-5.0472&-3.6800&-8.0911\\
-0.3772&-2.1179&-2.2113\\
-1.1831&-2.5827&0.2628
\end{array}\right]
$$
and
$$A_{2}=\left[\begin{array}{ccc}
2.5443&-0.5727&-6.9565\\
-0.4443&2.5444&-5.3728\\
0.8497&-1.5663&-0.5512
\end{array}\right].
$$
The time delays $\tau_1=2$ and $\tau_2=3.$
\par
Choose $T=200$ and $h=0.0001.$ Using the scheme, Eqs. (\ref{eq:state1-d}) to (\ref{eq:state4-d}), we have that
$$
\sum_{n=1}^{N}X_{n}h=\left[\begin{array}{ccc}
-0.2138&0.3819&0.0332\\
-0.2034&0.2524&0.2748\\
0.1296&-0.1438&-0.0613
\end{array}\right]
$$
and
$$
-(A_{0}+\sum_{j=1}^{m}A_{j})^{-1}=
\left[\begin{array}{ccc}
-0.2141&0.3813&0.0375\\
-0.2029&0.2500&0.2791\\
0.1287&-0.1423&-0.0612
\end{array}\right],
$$
respectively. They approach sufficiently. Thus Theorem \ref{th:P1} holds. System (\ref{eq:dde}) is asymptotically stable.
\par
Further, we obtain that $\sum_{n=1}^{N}||X_{n}||_{F} h =4.9285\quad \mbox{and}\quad \sum_{n=1}^{N}{||X_{n}||_{F}}^{2} h=14.3832.$
Thus Theorems \ref{th:P2} and \ref{th:P3} hold. The Frobenius norm of the fundamental matrix is given in Figure \ref{fig:Ex1}.
\begin{figure}
\centering
\includegraphics[width=8cm]{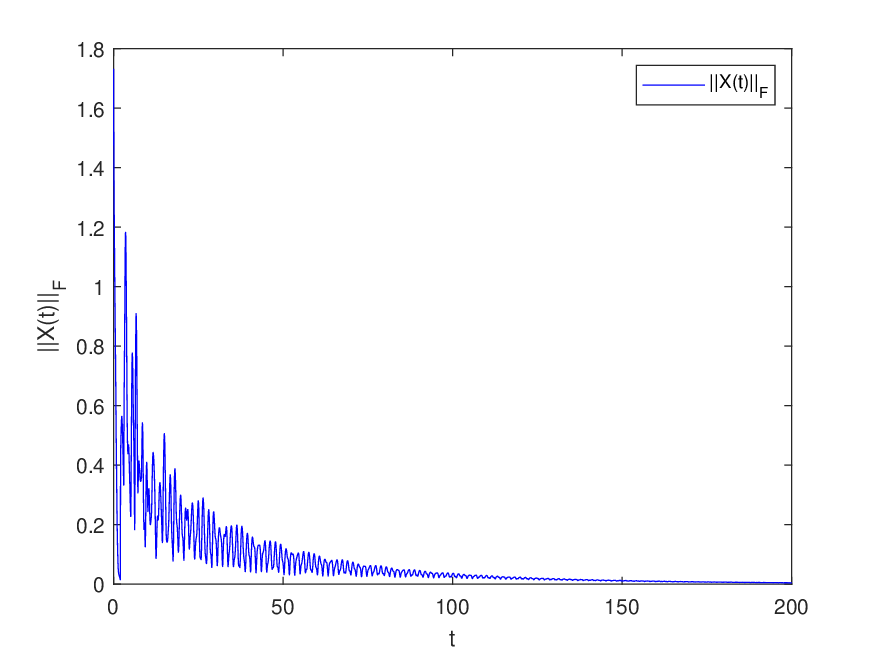}
\caption{$||X(t)||_{F}$: Frobenius Norm of Fundamental Matrix $X(t)$}
\label{fig:Ex1}
\end{figure}
\end{examp}
\par
\begin{examp}
For system (\ref{eq:dde}), let $m=2.$ The system matrices are
$$
A_{0}=\left[\begin{array}{cccc}
0.1535&0.0588&0.0417&0.1501\\
0&0&0.2000&-0.2000\\
-0.0802& -0.8381& -0.6125&0.2248\\
0.0267&0.6794&0.4708&-0.4749
\end{array}\right],
A_{1}=\left[\begin{array}{cccc}
-0.3538&0.5344&0.3015&-0.4100\\
0&-0.2000& -0.2000&0.2000\\
0.1307&-0.2515& -0.0522& -0.0350\\
0.0231&-0.4828& -0.5493&0.0450
\end{array}\right]
$$
and
$$
A_{2}=\left[\begin{array}{cccc}
0.1000 &-0.1000&0&0.2000\\
0&-0.2000& 0.1000&0\\
0&0.1000&0&0.300\\
0&-0.3000& 0.1000&0
\end{array}\right].
$$
The time delays $\tau_1=2$ and $\tau_2=3.$
\par
\par
Choose $T=100$ and $h=0.001.$ Using the scheme, Eqs. (\ref{eq:state1-d}) to (\ref{eq:state4-d}), we have that
$$
\sum_{n=1}^{N}X_{n}h=\left[\begin{array}{cccc}
18.9462&-2.6632&9.6488&8.3471\\
0.5533&1.6666&0.5544&0.5543\\
2.2143&-3.3299&2.2170&2.2177\\
2.1727&-0.8762&1.0954&3.2698
\end{array}\right]
$$
and
$$
-(A_{0}+\sum_{j=1}^{m}A_{j})^{-1}=
\left[\begin{array}{cccc}
19.0055&-2.6824&9.6805&8.3812\\
0.5556&1.6669& 0.5556& 0.5556\\
2.2225&-3.3326&2.2225& 2.2225\\
2.1791& -0.8783& 1.0989& 3.2745
\end{array}\right],
$$
respectively. They approach sufficiently. Thus Theorem \ref{th:P1} holds. System (\ref{eq:dde}) is asymptotically stable.
\par
Further, we obtain that $\sum_{n=1}^{N}||X_{n}||_{F} h =27.0951 \quad \mbox{and}\quad \sum_{n=1}^{N}{||X_{n}||_{F}}^{2} h=28.4794.$
Thus Theorems \ref{th:P2} and \ref{th:P3} hold. The Frobenius norm
of the fundamental matrix is given in Figure \ref{fig:Ex2}.
\begin{figure}
\centering
\includegraphics[width=8cm]{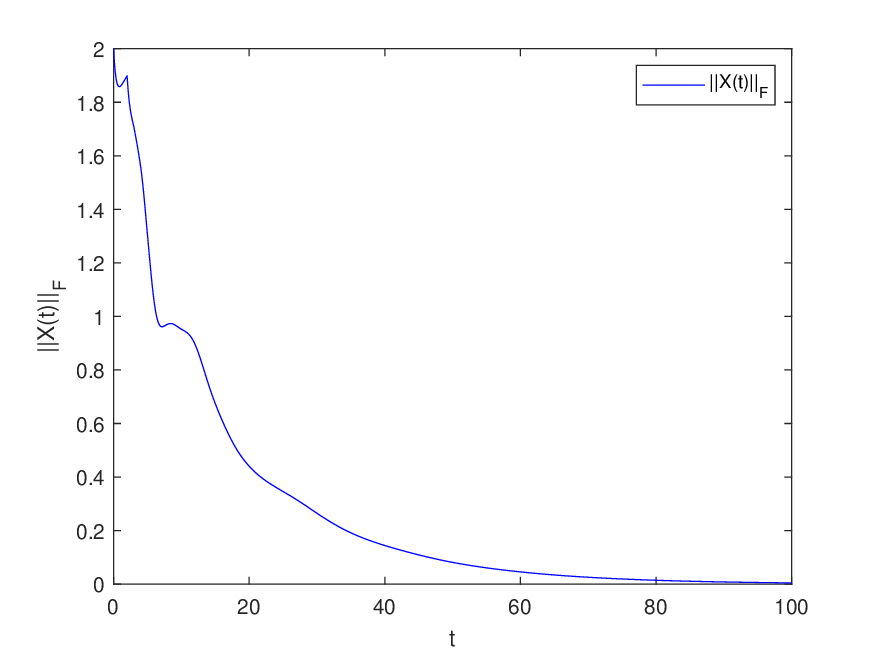}
\caption{$||X(t)||_{F}$: Frobenius Norm of Fundamental Matrix $X(t)$}
\label{fig:Ex2}
\end{figure}
\end{examp}
\section{Conclusions}
Based on integrals of the fundamental matrix, stability criteria are derived
which are necessary and sufficient conditions for stability of linear delay systems.
We prove that the necessary condition for stability of linear delay systems
in literature is also sufficient. Finally, numerical examples are given to
illustrate the main results.
\section{Acknowledgements}
This work is supported by the National Natural Science Foundation
of China (12271340 and 12371399). 

\end{document}